\documentclass[11pt]{article}

\def\corresponds{{\lower.2ex\hbox{=}}{\rm\kern-.75em^\triangle}}
\def\succsim{\succ\kern-.9em_\sim\kern.3em}
\def\precsim{\prec\kern-1em_\sim\kern.3em}
\def\slantfrac#1#2{\kern1em^{#1}\kern-.3em/\kern-.1em_{#2}}
\def\lfrac#1#2{{}^{#1\!}\kern-.0em/_{#2}}

\usepackage{amsfonts}

\begin{document}

\begin{center}
\begin{tabular}{c}
\hline
\rule[-5mm]{0mm}{15mm} {\Huge \sf Convergence Acceleration Techniques} \\
\hline
\end{tabular}
\end{center}
\vspace{0.2cm}
\begin{center}
U. D. Jentschura$^{1)}$, S. V. Aksenov$^{2)}$, P. J. Mohr$^{3)}$,
M. A. Savageau$^{2)}$, and G. Soff$^{1)}$
\end{center}
\vspace{0.2cm}
\begin{center}
$^{1)}$ {\em Institute of Theoretical Physics, \\
Dresden University of Technology, 01062 Dresden, Germany} \\[2ex]
$^{2)}$ {\em Department of Microbiology and Immunology, \\
University of Michigan, Ann Arbor, MI 48109, USA} \\[2ex]
$^{3)}$ {\em National Institute of Standards and Technology, \\
Mail Stop 8401, Gaithersburg, MD 20899, USA}
\end{center}
\vspace{0.3cm}
\begin{center}
\begin{minipage}{11.8cm}
{\underline{Abstract}}
This work describes numerical methods that are useful in many areas:
examples include statistical modelling (bioinformatics, computational
biology), theoretical physics, and even pure mathematics. The methods are
primarily useful for the acceleration of slowly convergent and the
summation of divergent series that are ubiquitous in relevant
applications. The computing time is reduced in many cases by orders of
magnitude.
\end{minipage}
\end{center}

\vspace{0.6cm}

Numerical approximation and analysis\\

\newpage

\section{CONVERGENCE ACCELERATION\\ (A BRIEF OVERVIEW)}

What does it mean to ``accelerate convergence''? The answer is that one
can do better than adding an infinite series term by term if the goal is
to get its sum to some specified numerical accuracy. This may appear as a
paradox, but the truth is that several powerful techniques have been
developed to that end since the arrival of the computer. The secret is to
use hidden information in trailing digits of partial sums of the input
series, to make assumptions on the form of the truncation error, and to
subsequently eliminate that error by a suitably chosen algorithm. Success
is judged by numerical experiments, and performance is the target.

A rather famous example for a problematic slowly convergent series is  
the Dirichlet series for the Riemann zeta function of argument 2,
\begin{equation}
\zeta(2) = \sum_{n=1}^{\infty} \frac{1}{n^2}
\end{equation}
whose terms all have the same (positive) sign.  It can be shown that a
tenfold increase in the number of added terms improves the accuracy of the
total sum only by a single decimal. When applying ``usual'' convergence
acceleration methods like the epsilon algorithm~\cite{Wy1956a} to {\em
nonalternating} series, severe numerical instabilities are more likely the
rule than the exception, and in general, more sophisticated algorithms
have to be sought. We have found that the {\em combination} of two
transformations leads to convincing numerical results in many
applications. The two steps are: (i) a Van Wijngaarden transformation
which transforms the nonalternating input series into an alternating
series, and (ii) the acceleration of the Van Wijngaarden transformed
series by a delta transformation. Details of the two transformations can be
found in~\cite{JeMoSoWe1999}, and further developments will
be published in~\cite{AkSaJeBeSoMo2003}.  
Here, we are just going to state that the partial sums of the
input, denoted $s_n$, are double-transformed to a series of transforms
$\delta_n^{(0)}(1,{\mathbf S}_n)$ according to the 
\begin{center}
COMBINED NONLINEAR--CONDENSATION TRANSFORMATION (CNCT), 
\end{center} 
that is $s_n \to \delta_n^{(0)}(1,{\mathbf S}_n)$, where the
$\delta_n^{(0)}(1,{\mathbf S}_n)$ exhibit much faster convergence than the
input data $s_n$. The example in Table 1 concerns the evaluation of ${\rm
Li}_3 (0.99999)$ to a relative accuracy of $10^{-15}$ already in twelfth
transformation order.  This result can also be achieved by term-by-term
summation of the defining series of ${\rm Li}_3$ -- in this case, however,
about 100 million terms are required.

\begin{center}
\begin{minipage}{7.5cm}
\begin{center}
{\sf Table 1.} Evaluation of 
$10^{-1}\,{\rm Li}_3 (0.99999)$ with the CNC
transformation~\cite{JeMoSoWe1999}.\\[0.5ex]
\begin{tabular}{lr}%
\hline
\hline
$n$ & \multicolumn{1}{l}{${\delta}_{n}^{(0)} \bigl(1, {\bf S}_0\bigr)$} \\
\hline%
 0 & \underline{0.1}33~331~333~415~539 \\
 1 & \underline{0.120}~474~532~168~000 \\
 2 & \underline{0.120~1}76~326~936~846 \\
 3 & \underline{0.120~204}~748~497~388 \\
 4 & \underline{0.120~204~0}79~128~106 \\
 5 & \underline{0.120~204~045}~387~208 \\
 6 & \underline{0.120~204~045}~378~284 \\
 7 & \underline{0.120~204~045~43}4~802 \\
 8 & \underline{0.120~204~045~438}~553 \\
 9 & \underline{0.120~204~045~438~7}26 \\
10 & \underline{0.120~204~045~438~733} \\
11 & \underline{0.120~204~045~438~733} \\
12 & \underline{0.120~204~045~438~733} \\
\hline%
exact & \underline{0.120~204~045~438~733} \\
\hline
\hline
\end{tabular}
\end{center}
\end{minipage}
\end{center}

\section{APPLICATIONS IN BIOPHYSICS}

The theoretical description of biological processes is unthinkable today
without extensive statistical analysis, and concurrently, fields like
``bioinformatics'' and ``computational biology'' are emerging. Several
important mathematical functions needed in the theory of statistical
distributions are represented by slowly convergent series, and their
computation can benefit to a large extent from using the CNCT proposed
here. Examples include the {\em discrete Zipf-related distributions} whose
probability mass functions are represented by the terms of infinite series
defining Riemann zeta, generalized zeta, and polylogarithm functions and
whose total probability is calculated with these functions. These
distributions are used in statistical analysis of biological sequences (of
RNA, DNA and protein molecules) and occurence analysis of folds of
proteins~\cite{AkSaJeBeSoMo2003}. The generalized representation of these
distributions was shown to be in the form of the Lerch distributional
family~\cite{AkSaJeBeSoMo2003,ZoAl1995,AkSa2001}, which requires
calculation of Lerch's $\Phi$ transcendent. The $\Phi$ transcendent 
is given by the following power series,
\begin{equation}
\label{lerch}
\Phi(z,s,v) = \sum_{n=0}^\infty \frac{z^n}{(n+v)^s}\,.
\end{equation}
For $|z| < 1$, $|z| \approx 1$, the power series is very slowly
convergent. Of particular importance is the case of a real
argument $x \equiv z$. In the region $x \approx 1$, 
we found that the application of the 
CNC transformation~\cite{JeMoSoWe1999} leads to a significant 
acceleration of the convergence, whereas for $x \approx -1$,
numerical problems can be solved by the direct application of the 
delta transformation [see Eq.~(8.4-4) of~\cite{We1989}]
to the defining series~(\ref{lerch}). The Riemann zeta, generalized zeta,
and polylogarithm functions are special cases of Lerch's transcendent.
Further applications of this special function include the quantile
function of {\em continuous {\em S} distributions}~\cite{HeSo2001}.
Finally, the evaluation of several hypergeometric functions and related
hypergeometric distributions can be significantly enhanced using the CNCT.
Needless to say, a fast and accurate computation of these special
functions is of crucial importance in calculating various basic properties
of these distributions, including moments, cumulative distribution
functions and quantiles, and parameter estimations.

\section{APPLICATIONS IN THEORETICAL PHYSICS}

We will briefly mention that various 
{\em long-standing problems in theoretical
physics} have recently been solved using computational methods based on
convergence acceleration techniques. Examples include quantum
electrodynamic bound state calculations~\cite{JeMoSo1999}, which yield a
theoretical description of the most accurate physical measurements today
(in some cases, laser spectroscopy has reached a relative accuracy of
$10^{-14}$~\cite{NiEtAl2000}). Therefore, the calculations are of
importance for the test of fundamental quantum theories and for the
determination of fundamental physical constants. Further applications
include the evaluation of quantum corrections to Maxwell equations, which
are given by the ``quantum electrodynamic effective
action''~\cite{JeGiVaLaWe2001}. This object is representable by a slowly
convergent series and is phenomenologically important in the description of
various astrophysical processes.

We also report that it is possible, in combining analytic results
obtained in~\cite{JeSoMo1997} with numerical techniques based on the 
CNCT, to evaluate the so-called Bethe logarithm in hydrogen to 
essentially arbitrary precision. Specifically, we obtain -- for the 
4P state -- the result
\begin{equation}
\ln k_0(4 {\mathrm P}) = -0.041~954~894~598~085~548~671~037(1)
\end{equation}
which is 9 orders of magnitude more accurate than the 
latest and most precise calculation recorded so far in the
literature~\cite{DrSw1990}.

\section{APPLICATIONS IN MATHEMATICS}

As far as mathematics is concerned, we will quote from~\cite{Ba1994tech}:
``In April 1993, Enrico Au-Yeung, an undergraduate at the University of
Waterloo, brought to the attention of [David Bailey's] colleague Jonathan
Borwein the curious fact that
\begin{eqnarray*}
\label{pisum}
\sum_{k=1}^\infty \left(1 + {1 \over 2} + \cdots + {1 \over k}\right)^2
  \, k^{-2} &=& \frac{17 \pi^4}{360}
\end{eqnarray*}
based on a computation to 500,000 terms.  Borwein's reaction was to
compute the value of this constant to a higher level of precision in order
to dispel this conjecture.  Surprisingly, his computation to 30 digits
affirmed it. [David Bailey] then computed this constant to 100 decimal
digits, and the above equality was still affirmed.''

Many formulas similar to (\ref{pisum}) have subsequently been established
by rigorous proof~\cite{BaBoGi1994}. Using the CNCT, it is easy to
calculate the sum (\ref{pisum}) to 200 digits, based on multiprecision
arithmetic~\cite{Ba1994tech} and a Linux personal computer, within a few
hours. In calculating the specific case (\ref{pisum}) to an accuracy of
200 decimals, which {\em a priori} requires the calculation of about
$10^{205}$ terms of the series, we report that roughly 84~000 terms are
sufficient when employing the CNCT. This corresponds to an 
acceleration of
the convergence by roughly 200 orders of magnitude.

\section{CONCLUSION}

The CNCT has become useful in a wide variety of application areas which
extend beyond the original scope of the
transformation~\cite{JeMoSo1999}.
Details of the implementation of the algorithm, in the three
languages C, Fortran and Mathematica~\cite{disclaimer} 
will be presented at the 
conference. Sample files will also be made available for
internet download at~\cite{AkJeURL}. We have recently
investigated further potential applications of the
algorithms described here, such as 
the evaluation of generalized hypergeometric functions
which can be of exquisite practical importance,
with rather promising results. The rather general 
applicability of the convergence accelaration methods
makes them very attractive tools in scientific computing.


\begin{thebibliography}{99}
\bibitem{disclaimer}
Certain commercial equipment, instruments, or materials are identified in this
paper to foster understanding.
Such identification does not imply recommendation or endorsement by the
National Institute of Standards and Technology, nor does it imply that
the materials or equipment identified are necessarily the best
available for the purpose.
\end{thebibliography}
\end{document}